\documentclass{amsart}
\usepackage{graphicx}
\usepackage{amssymb}

\newtheorem{thm}{Theorem}[section]

\newtheorem{prop}[thm]{Proposition}
\theoremstyle{definition}

\theoremstyle{remark}
\newtheorem{rem}[thm]{Remark}
\theoremstyle{conjecture}

\theoremstyle{example}
\newtheorem{ex}[thm]{Example}

\begin{document}

\title[Singular structures]{Singular structures in solutions to the Monge-Amp\`{e}re equation with point masses}

\dedicatory{Dedicated with admiration to Professor Neil Trudinger on the occasion of
  his 80th birthday.}

\author{Connor Mooney}
\address{Department of Mathematics, UC Irvine}
\email{\tt mooneycr@math.uci.edu}
\author{Arghya Rakshit}
\address{Department of Mathematics, UC Irvine}
\email{\tt arakshit@uci.edu}

\begin{abstract}
We construct new examples of Monge-Amp\`{e}re metrics with polyhedral singular structures, motivated by problems related to the optimal transport of point masses and to mirror symmetry. We also analyze the stability of the singular structures under small perturbations of the data given in the problem under consideration.
\end{abstract}
\maketitle

\section{Introduction}
Monge-Amp\`{e}re metrics with singularities appear in a variety of contexts, including mirror symmetry (in connection with the Strominger-Yau-Zaslow conjecture, see e.g. \cite{CL2}, \cite{JX}, \cite{Li}, \cite{Lo}, \cite{LYZ}) and in the optimal transport of singular measures. By a Monge-Amp\`{e}re metric we mean the Hessian of a convex solution to $\det D^2u = 1$. In \cite{M2} a robust method was developed to construct examples of such metrics with Y-shaped and polyhedral singular structures in three and four dimensions, based on solving a certain obstacle problem. The extension of the main  result in \cite{M2} to higher dimensions, stated as a conjecture (Conjecture 1.4 in that paper), was hindered by the lack of a well-developed regularity theory for the obstacle problem considered in that work. The purpose of this paper is to prove Conjecture 1.4 from \cite{M2} using a simplified approach which avoids the use of delicate free boundary regularity results, to analyze the stability of the singular structures appearing in these examples, and to suggest further research directions motivated by the connections of our examples to the aforementioned areas.

\vspace{2mm}

Our main result is:

\begin{thm}\label{Main}
Let $P \subset \mathbb{R}^n$ be a compact convex polytope, and let $\Gamma_k$ denote its $k$-skeleton. Then there exists a convex function $u : \mathbb{R}^n \rightarrow \mathbb{R}$ such that
$$\Gamma_{\left\lceil \frac{n}{2} - 1\right\rceil} \subset \{u = 0\}, \quad u \in C^{\infty}\left(\mathbb{R}^n \backslash\Gamma_{\left\lceil \frac{n}{2} - 1\right\rceil}\right), \quad \text{and} \quad \det D^2u = 1 + \sum_{q \in \Gamma_0} a_q\delta_q$$
for some coefficients $a_q > 0$.
\end{thm}

\noindent In particular, $u$ is singular on $\Gamma_{\left\lceil \frac{n}{2} - 1\right\rceil}$, and $\nabla u$ is in fact discontinuous there (see Remark \ref{GradJump}). Here $\lceil t \rceil$ denotes the smallest integer greater than or equal to $t$. In dimension $n = 2$ solutions to $\det D^2u = 1$ are locally strictly convex and smooth \cite{A}, so the examples proving Theorem \ref{Main} exhibit interesting singular structures away from the vertex set $\Gamma_0$ in dimensions three and larger.

In \cite{M2}, the approach to Theorem \ref{Main} (which was successful in dimensions four and smaller) was based on solving an obstacle problem by lowering super-solutions to the Monge-Amp\`{e}re equation $\det D^2u = 1$ while constraining them to lie above a polyhedral graph (the obstacle). In this paper we instead consider a ``dual" obstacle problem, where we raise sub-solutions 
to the equation from below while constraining them to lie below prescribed values at the vertices of $P$. In this way we can avoid using delicate regularity results from \cite{S} which were necessary for carrying out the previous approach.

\vspace{2mm}

We then study the stability of the singular structures in the solutions from Theorem \ref{Main} from two perspectives. First, global solutions on $\mathbb{R}^n$ to equations of the form
\begin{equation}\label{MassEqn}
\det D^2u = 1 + \sum_{i = 1}^M a_i \delta_{p_i}
\end{equation}
are asymptotic to quadratic polynomials \cite{CL}. Modulo affine invariance, the space of solutions to (\ref{MassEqn}) can be identified with an explicit orbifold parametrized by the mass sizes $a_i$ and the mass locations $p_i$ (see \cite{JX}).
It is natural to ask about the geometry and topology of the set in this moduli space which corresponds to solutions $u$ that are singular away from $\{p_i\}_{i = 1}^M$. Our proof of Theorem \ref{Main} shows that
this set is not small. In particular, it has nonempty interior:

\begin{thm}\label{Stability1}
Let $u$ be one of the examples constructed in the proof of Theorem \ref{Main}, and assume that it solves
$$\det D^2u = 1 + \sum_{i = 1}^M a_i \delta_{p_i}$$
for some $a_i > 0$ and $p_i \in \mathbb{R}^n$. If $\sum_{i = 1}^M (|\tilde{a}_i - a_i| + |\tilde{p}_i - p_i|)$ is sufficiently small, then the global solution $\tilde{u}$ to 
$$\det D^2\tilde{u} = 1 + \sum_{i = 1}^M \tilde{a}_i \delta_{\tilde{p}_i}$$
which is asymptotic to the same quadratic polynomial as $u$ is singular on the faces of the polytope with vertices $\{\tilde{p_i}\}_{i = 1}^M$ that have dimension smaller than $n/2$, and $\tilde{u}$ is smooth elsewhere.
\end{thm}

Second, the Legendre transform $u^*$ of one of the examples $u$ constructed in Theorem \ref{Main} can be viewed as the potential of the optimal transport map (with quadratic cost) which pushes forward the Lebesgue measure $dx$ in a bounded domain $\Omega^* \subset \mathbb{R}^n$ to the measure 
$$\nu = dx + \sum_{q \in \Gamma_0} a_q\delta_q$$ 
on $\Omega := \nabla u^*(\Omega^*)$. The dual optimal transport plan $\nabla u$ is discontinuous on $\Gamma_{\left\lceil \frac{n}{2} - 1\right\rceil}$, despite $\nu$ being regular away from $\Gamma_0$. Since $\Omega^*$ can be taken convex, the singularities are a result of the presence of Dirac masses in $\nu$ and not of the boundary geometry of $\Omega^*$ (if $\nu$ had a smooth positive density and $\Omega^*$ were convex, then the optimal transport map from $\nu$ in $\Omega$ to $dx$ in $\Omega^*$ would be smooth \cite{COT}). It is natural to ask if the discontinuities persist when the data of the problem (the measures) are perturbed. Our last result shows that they do:

\begin{thm}\label{Stability2}
Let $u$ be one of the examples constructed in the proof of Theorem \ref{Main}, and assume that it solves
$$\det D^2u = dx + \sum_{i = 1}^M a_i \delta_{p_i} := \nu$$
for some $a_i > 0$ and $p_i \in \mathbb{R}^n$.
Let $\Omega^*$ be a bounded convex domain containing $\nabla u(P)$, and let $\Omega = \nabla u^*(\Omega^*)$. If $\sum_{i = 1}^M (|\tilde{a}_i - a_i| + |\tilde{p}_i - p_i|)$ is sufficiently small and 
$$\tilde{\nu} := dx + \sum_{i = 1}^M \tilde{a}_i \delta_{\tilde{p}_i}$$
satisfies the mass balance condition $\tilde{\nu}(\Omega) = \nu(\Omega)$ (that is, $\sum_{i = 1}^M \tilde{a}_i = \sum_{i = 1}^M a_i$), then the Legendre transform $\tilde{u}$ of the potential $\tilde{u}^*$ of the optimal transport from the Lebesgue measure in $\Omega^*$ to $\tilde{\nu}$ in $\Omega$ satisfies that $\nabla \tilde{u}$ is discontinuous 
on the faces of the polytope with vertices $\{\tilde{p_i}\}_{i = 1}^M$ that have dimension smaller than $n/2$, and $\tilde{u}$ is smooth elsewhere in $\Omega$.
\end{thm}

The intuition for Theorems \ref{Stability1} and \ref{Stability2} is that if the mass locations $p_i$ are close to one another and the masses $a_i$ are large, then the masses ``communicate" and singularities are generated between them in optimal transport maps (in dimensions three and higher, at least). If on the other hand the masses are far from one another in comparison to the mass sizes, they do not communicate and the transport maps are smooth away from the masses (see Example \ref{SmoothBetween}).

\vspace{2mm}

The paper is organized as follows. In Section \ref{Preliminaries} we recall the notion of Monge-Amp\`{e}re measure, solve an obstacle problem, and recall a family of useful Pogorelov-type singular solutions. In Section \ref{ProofMain} we prove Theorem \ref{Main}. In Section \ref{Stability} we prove Theorems \ref{Stability1} and \ref{Stability2}. Finally, in Section \ref{OpenProbs} we list and discuss some open questions motivated by this work.

\section*{Acknowledgments}
The authors gratefully acknowledge the support of NSF grant DMS-1854788, NSF CAREER grant DMS-2143668, an Alfred P. Sloan Research Fellowship, and a UC Irvine Chancellor's fellowship.


\section{Preliminaries}\label{Preliminaries}

In this section we recall the notion of Monge-Amp\`{e}re measure, solve an obstacle problem for the Monge-Amp\`{e}re equation, define a family of Pogorelov-type singular solutions to the Monge-Amp\`{e}re equation,
and recall a regularity result from \cite{C} which bounds the dimension of a singularity in a solution to the Monge-Amp\`{e}re equation.

\subsection{Monge-Amp\`{e}re Measure}
To a convex function $v$ on a domain $\Omega \subset \mathbb{R}^n$ we associate a Borel measure $Mv$ on $\Omega$, called the Monge-Amp\`{e}re measure of $v$. It satisfies
$$Mv(E) = |\partial v(E)|$$
for any Borel set $E \subset \Omega$, where $\partial v$ denotes the subgradient of $v$. When $v \in C^2$ we have $Mv = \det D^2v \,dx.$
Given a Borel measure $\mu$ on $\Omega$, we say that $v$ is an Alexandrov solution to the Monge-Amp\`{e}re equation $\det D^2v = \mu$
if  $Mv = \mu.$ 

Alexandrov solutions are closed under uniform convergence: if convex functions $v_k$ converge locally uniformly in $\Omega$ to
$v$, then their Monge-Amp\`{e}re measures $Mv_k$ converge weakly to $Mv$. 

Finally, given a bounded convex domain $\Omega \subset \mathbb{R}^n$ and a finite Borel measure $\mu$ on $\Omega$, the Dirichlet problem
$$\begin{cases}
Mv = \mu \text{ in } \Omega, \\
v|_{\partial \Omega} = \varphi
\end{cases}$$
is solvable in $C\left(\overline{\Omega}\right)$ provided e.g. $\varphi$ is linear, or $\Omega$ is strictly convex and $\varphi$ is continuous.
For proofs of these results see \cite{Gut}.


\subsection{Obstacle Problem}\label{ObstacleProblem}

We now solve an obstacle problem. The data are a bounded strictly convex domain $U \subset \mathbb{R}^n$, boundary data $\varphi \in C\left(\overline{U}\right)$, an obstacle $g : \overline{U} \rightarrow \mathbb{R} \cup \{+ \infty\}$ which is lower semicontinuous and satisfies $g > \varphi$ on $\partial U$, and a finite Borel measure $\mu$ on $U$. We define the class of functions $\mathcal{F}$ by
$$\mathcal{F} := \left\{v : v \in C\left(\overline{U}\right) \text{ convex},\, v \leq g \text{ in } U,\, v|_{\partial U} = \varphi,\, Mv \geq \mu\right\}.$$
We show:
\begin{prop}\label{OP}
The set $\mathcal{F}$ is non-empty, the function
$$u := \sup_{\mathcal{F}} v$$ 
is in $\mathcal{F}$, and 
$$Mu = \mu \text{ in } \{u < g\} \cap U.$$
\end{prop}

\vspace{2mm}

\begin{proof}[{\bf Proof of Proposition \ref{OP}}]
Let $u_0$ be the solution in $C\left(\overline{U}\right)$ to
$$\begin{cases}
Mu_0 = \mu \text{ in } U, \\
u_0|_{\partial U} = 0,
\end{cases}$$
and let $\varphi_0$ be the convex envelope of $\min\{\varphi,\,g\}$. Then $u_1 := u_0 + \varphi_0 \in \mathcal{F}$. Let $\varphi_1$ be the convex envelope of the boundary data of $\varphi$ (the supremem of affine functions that are lower than $\varphi$ on $\partial U$), which satisfies $\varphi_1 \in C\left(\overline{U}\right),\, \varphi_1|_{\partial U} = \varphi,$ and $M\varphi_1 = 0$ in $U$ (see e.g. \cite{Gut}). Using that $\mathcal{F}$
is closed under taking maxima and under uniform convergence, it is not hard to construct an increasing sequence of functions $u_k \in \mathcal{F}$ which satisfy
$$u_1 \leq u_k \leq \varphi_1$$
for all $k$ and tend uniformly to $u \in \mathcal{F}$. 

To conclude we show that for any $x$ in the open set $\{u < g\}$, there exists $r_x > 0$ such that $Mu = \mu$ in $B_r(x)$ for all $r < r_x$. Since such balls generate the Borel $\sigma$-algebra, this will complete the proof. Let $w$ be the solution to
$$\begin{cases}
Mw = \mu \text{ in } B_r(x), \\
w|_{\partial B_r(x)} = u.
\end{cases}
$$
For $r$ small we have $w < g$ in $B_r(x)$, and from the maximum principle we have $u \leq w$. Replacing $u$ by $w$ in $B_r(x)$ we obtain a function in $\mathcal{F}$, hence $u = w$ in $B_r(x)$ and we are done.
\end{proof}

\begin{rem}\label{ModifiedPerron}
We can also write $u$ as the supremum of functions in 
$$\mathcal{\tilde{F}} := \left\{\tilde{v} : \tilde{v} \in C\left(\overline{U}\right) \text{ convex},\, \tilde{v} \leq g \text{ in } U,\, \tilde{v}|_{\partial U} \leq \varphi,\, M\tilde{v} \geq \mu\right\}.$$
Indeed, for any function $\tilde{v} \in \mathcal{\tilde{F}}$ there is a function $v \in \mathcal{F}$ such that $v \geq \tilde{v}$, given by the maximum between $\tilde{v}$ and the function $u_1 \in \mathcal{F}$ defined in the proof of Proposition \ref{OP}.
\end{rem}


\subsection{Barriers}\label{Barriers}
We now define a useful family of Pogorelov-type barriers constructed in \cite{C}.
We denote points in $\mathbb{R}^n$ by $(x,\,y)$ with $x \in \mathbb{R}^{n-k}$ and $y \in \mathbb{R}^k$. For $n \geq 3$ and $1 \leq k < \frac{n}{2}$, define the function $w_{n,\,k}$ on $\mathbb{R}^n$ by
\begin{equation}\label{PogoEx}
w_{n,\,k}(x,\,y) = C(n) |x|^{2-2k/n}(1+|y|^2).
\end{equation}
For $C(n)$ sufficiently large we have
$$\det D^2w_{n,\,k} \geq 1$$
(in the Alexandrov sense) in the slab $\{|y| < \rho_n\}$ for some $\rho_n > 0$. We omit the calculation, which is straightforward using coordinates that are polar in $x$ and $y$. 

\vspace{2mm}

For $n \geq 1$ we also let
\begin{equation}\label{OneDelta}
W_n(x) := \int_0^{|x|} \left(1 + s^n\right)^{\frac{1}{n}}\,ds,
\end{equation}
which solves
\begin{equation}\label{FlatSuper}
\det D^2W_n = 1 + |B_1|\delta_0
\end{equation}
in the Alexandrov sense. It also satisfies
\begin{equation}\label{GrowthRadial}
W_n(x) - \frac{1}{2}|x|^2 = \begin{cases}
O(|x|), \quad n = 1 \\
O(|\log|x||), \quad n = 2 \\
c(n) + O(|x|^{2-n}), \quad n \geq 3
\end{cases}
\end{equation}
for some constants $c(n) > 0$, and 
\begin{equation}\label{DistGrowth}
W_n(x) \geq |x|
\end{equation}
on $\mathbb{R}^n$.

\subsection{A regularity result}\label{Regularity}
To conclude the section we recall a useful bound on the dimension of a singularity appearing in a solution to the Monge-Amp\`{e}re equation (\cite{C}, see also \cite{M1} for a short proof).

\begin{prop}\label{MaxDimension}
Assume that $\det D^2u \geq 1$ in the Alexandrov sense in a domain $U \subset \mathbb{R}^n$, and let $L$ be a supporting linear function to $u$. Then
$$\text{dim}\{u = L\} < \frac{n}{2}.$$
\end{prop}

\noindent The examples $w_{n,\,k}$ show that this bound is optimal.

\section{Proof of Theorem \ref{Main}}\label{ProofMain}
Below we will use the following observation: there exists some $\delta > 0$ depending on $P$ such that, for any face $F \in \Gamma_k$ with $k < n$, there is an affine function $L$ that satisfies 
$$|\nabla L| = 1,\, L|_F = 0,\, \text{ and } L \leq -\delta \text{ on } \Gamma_0 \backslash F.$$
We also assume that $n \geq 3$, in view of the local regularity theory for the Monge-Amp\`{e}re equation in two dimensions mentioned in the introduction.

\begin{proof}[{\bf Proof of Theorem \ref{Main}}]
After a translation we may assume that $0 \in P$. By quadratic rescaling we may replace $P$ with $\epsilon_0P$ for $\epsilon_0 > 0$ small depending on $n,\,P$ to be chosen. Let $u_R$ be the solution to the obstacle problem from Section \ref{ObstacleProblem} with 
$$U = B_R,\, \varphi = W_n + 1,\, \mu = dx,$$ 
and 
$$g(x) = \begin{cases}
0,\, x \in \epsilon_0\Gamma_0 \\
+ \infty,\, \text{otherwise}.
\end{cases}
$$
By the maximum principle and the fact that $u_R(0) \leq 0$ we have $u_R \leq \varphi$. Here and below, we will let $C$ denote a large constant depending on $n$ and $P$. By the definition of $u_R$, provided $C$ is chosen sufficiently large we have that 
$$\epsilon_0 P \subset \{W_n - C\epsilon_0 < 0\},$$ hence 
$$W_n - C\epsilon_0 \leq u_R$$ 
in $B_R$ for all $R$ (see Remark \ref{ModifiedPerron}).

For $k < \frac{n}{2}$ and any face $F \in \Gamma_k$, choose an affine function $L$ such that 
$$|\nabla L| = 1,\, L|_{\epsilon_0F} = 0,\, \text{ and } L < -\delta\epsilon_0$$ 
at all points in $\epsilon_0\Gamma_0 \backslash \epsilon_0 F$. Let $z_0 \in \epsilon_0F$. For some rotation $O$, the function 
$$B(x) := w_{n,\,k}(O(x -z_0))$$ 
vanishes on $\epsilon_0 F$. Fixing $\rho(n)$ small, we have that
$$B - L/2 < W_n - C\epsilon_0$$ 
on $\partial B_{\rho}$ provided $\epsilon_0$ is small, using that $W_n(x) \geq |x|$ and that
$$B(x) \leq C(n)|x-z_0|^{1+1/n}$$ 
for $|x - z_0| < 1$. Finally, taking $\epsilon_0$ smaller if necessary, we have 
$$B - L/2 \leq C\epsilon_0^{1+1/n} - \delta\epsilon_0/2 < 0$$ 
at all points in $\epsilon_0\Gamma_0 \backslash \epsilon_0 F$. We conclude that the function obtained by replacing $W_n - C\epsilon_0$ by $\max\{W_n - C\epsilon_0,\, B - L/2\}$ in $B_{\rho}$ is
in the class $\tilde{\mathcal{F}}$ defined in Remark \ref{ModifiedPerron}, hence $u_R = B - L/2 = 0$ on $\epsilon_0 F$.

Using that $W_n - C\epsilon_0 \leq u_R \leq W_n + 1$ for all $R$, we may take a sequence of radii $R_j$ tending to infinity such that the corresponding $u_{R_j}$ converge locally uniformly to a global convex function $u$
which solves $\det D^2u = 1$ away from $\epsilon_0\Gamma_0$ and vanishes on $\epsilon_0\Gamma_k$ for all $k < n/2$. We claim that $u$ is smooth otherwise. Outside the  polytope this follows from results in \cite{CL}, which say that $u$ is strictly convex (hence smooth) outside the convex hull of $\epsilon_0\Gamma_0$ (that is, $\epsilon_0P$). To finish, we claim that $\{u < 0\}$ contains the interiors of all faces of $\epsilon_0P$ of dimension $n/2$ or larger. Indeed, if $u$ vanishes at an interior point of such a face, then $u$ vanishes in the whole face by convexity, which contradicts Proposition \ref{MaxDimension}. Since $\det D^2u = 1$ in $\{u < 0\}$, the function $u$ is smooth in $\{u < 0\}$ by classical results (\cite{P}, \cite{CY}) and the proof is thus complete.
\end{proof}

\begin{rem}\label{GradJump}
It is in fact true that $\nabla u$ is discontinuous on $\Gamma_{\left\lceil \frac{n}{2} - 1\right\rceil}$. Indeed, in the proof of Theorem \ref{Main} we can replace $w_{n,\,k}$ by appropriate rescalings of different Pogorelov-type sub-solutions of the form
$$\tilde{w}_{n,\,k} = |x| + |x|^{\frac{n-k+1}{k+1}}(1+|y|^2),$$
which also vanish on the $k$-dimensional subspace $\{|x| = 0\}$ and have a Lipschitz singularity on this subspace.
\end{rem}

\section{Proofs of Theorems \ref{Stability1} and \ref{Stability2}}\label{Stability}

The barrier arguments in the proof of Theorem \ref{Main} show that the presence of singularities is robust under $C^0$ perturbations. By this we mean:

\begin{prop}\label{StabilityProp}
Let $u$ be an example constructed in the proof of Theorem \ref{Main}, and assume that $u$ solves
$$\det D^2u = 1 + \sum_{i = 1}^M a_i \delta_{p_i}.$$
If $\tilde{u}$ is a convex function defined in a neighborhood $N$ of $P$ such that $\det D^2\tilde{u} = 1$ away from points $\{\tilde{p}_i\}_{i = 1}^M$ with $\sum_{i = 1}^M |p_i - \tilde{p}_i|$ sufficiently small, and furthermore 
$\|u - \tilde{u}\|_{C^0(N)}$ is sufficiently small, then $\tilde{u}$ is singular on the faces of the polytope $\tilde{P}$ with vertices $\{\tilde{p}_i\}_{i = 1}^M$ of dimension smaller than $n/2$, and $\tilde{u}$ is smooth otherwise in a neighborhood $N' \subset N$ of $\tilde{P}$.
\end{prop}
\noindent Below we sketch the proof, suppressing the ``$\epsilon_0"$ from the proof of Theorem \ref{Main} for simplicity of notation.
\begin{proof}
By perturbing the barriers $B - L/2$ from the proof of Theorem \ref{Main} and applying the maximum principle, we see that $\tilde{u}$ is singular on each face $\tilde{F}$ of $\tilde{P}$ that has dimension smaller than $n/2$. More precisely, the convex envelope of the values that $\tilde{u}$ takes on the vertices of $\tilde{F}$ is linear on sub-regions that partition $\tilde{F}$ (see Example \ref{Subtlety} below). For each such sub-region, we can perturb $B-L/2$ by a small translation, rotation, and addition of an affine function to get a new barrier $\tilde{B} + \tilde{L}$ that agrees with the envelope on this sub-region, where $\tilde{B}$ is a rotation and translation of $w_{n,\,k}$ that is linear when restricted to the affine subspace containing $\tilde{F}$ and $\tilde{L}$ is linear and vanishes on $\tilde{F}$. The maximum principle and the convexity of $\tilde{u}$ imply that $\tilde{u}$ agrees with $\tilde{B} + \tilde{L}$ (in particular, is linear) on this sub-region. Key points are that $\tilde{B} + \tilde{L} \leq \tilde{u}$ at all vertices of $\tilde{F}$ by construction, and $\tilde{B} + \tilde{L}$ is close to $B - L/2$ which is less than $\tilde{u}$ at the remaining vertices of $\tilde{P}$ and on the boundary of the ``large" domain $N$.

As for regularity, assume another singularity happens in $\tilde{u}$. Its only extremal points can be some subset of the vertices $\{\tilde{p}_i\}_{i = 1}^M$, thus it is a polytope of dimension smaller than $n/2$ whose vertices are contained in $\{\tilde{p}_i\}_{i = 1}^M$. (Recall that a singularity has no extremal points on $N \backslash \cup_{i = 1}^M \tilde{p}_i$ \cite{C0}. To rule out the case of a singularity that extends from $N'$ to $\partial N$, use that $u$ is strictly convex outside of $P$ and that $\|\tilde{u} - u\|_{C^0(N)}$ is small). There is some $\mu > 0$ such that each such polytope which is not contained in a face of $\tilde{P}$ of dimension smaller than $n/2$ intersects $\{u < -\mu\}$, provided $\sum_{i = 1}^M |p_i - \tilde{p}_i|$ is small. Then $\tilde{u} < -\mu/2$ at such points provided $\|u - \tilde{u}\|_{C^0(N)}$ is small, giving a contradiction (on the singularity, $\tilde{u}$ is bounded between its values at the vertices which are small).
\end{proof}

\begin{proof}[{\bf Proofs of Theorems \ref{Stability1} and \ref{Stability2}}]
In view of Proposition \ref{StabilityProp}, it suffices to show that $\tilde{u}$ is close to $u$ in a neighborhood of $P$. Indeed, in the context of Theorem \ref{Stability1}, $\tilde{u}$ is smooth outside the convex hull of $\{\tilde{p}_i\}_{i = 1}^M$ by results in \cite{CL}, and in the context of Theorem \ref{Stability2} the function $\tilde{u}$ is smooth outside the convex hull of $\{\tilde{p}_i\}_{i = 1}^M$ by a small modification of the arguments in \cite{COT}.

To see this in the setting of Theorem \ref{Stability1}, assume that $u_k$ are the unique global solutions to 
$$\det D^2u_k = 1 + \sum_{i = 1}^M a^k_i\delta_{p^k_i}$$
that are asymptotic to the same quadratic polynomial as $u$ (after performing an affine change of variable and adding a linear function we may assume this is $|x|^2/2$), with $\sum_{i = 1}^M (|a^k_i - a_i| + |p^k_i - p_i|)$ tending to $0$ as $k \rightarrow \infty$. (See \cite{JX} for a discussion of the existence and uniqueness of solutions to this global problem). As shown in \cite{JX}, the functions 
$$v_k = \frac{1}{M} \sum_{i = 1}^M (\lambda_i^k)^2W_n((\cdot - p^k_i)/\lambda_i^k),$$ 
where 
$$\lambda^k_i = M(a_i^k/|B_1|)^{1/n},$$
are, up to adding quadratics with uniformly bounded (in $k$) coefficients, sub-solutions to the problem solved by $u_k$. These satisfy that $v_k \geq |x|^2/2 - K$ for some $K > 0$ and all $k$, and $|v_k - |x|^2/2| \leq K|x|^{2-n}$. By the maximum principle we have that $u_k \leq |x|^2/2$ for all $k$. We conclude from the inequality $v_k \leq u_k \leq |x|^2/2$ that any subsequence of $\{u_k\}$ has itself a locally uniformly convergent subsequence whose limit must (by uniqueness) be $u$.

In the optimal transport setting, one can use use the stability of optimal transport maps, applied to the Legendre transform. More precisely, let $a^k_i$ and $p^k_i$ be as above satisfying in addition the balance condition 
$$\sum_{i = 1}^M a^k_i = \sum_{i = 1}^M a_i,$$
and let $\nabla u_k^*$ be the optimal transport maps from the Lebesgue measure $dx$ in $\Omega^*$ to the measure
$$\nu_k = dx +  \sum_{i = 1}^M a^k_i\delta_{p^k_i}$$
in $\Omega$. Then the maps $\nabla u_k^*$ converge in measure to $\nabla u^*$ (see \cite{DP} Theorem 1.14), which along with their uniform boundedness in $\Omega^*$ implies the $C^0$ convergence of $u_k^*$ to $u^*$, up to adding constants to $u_k^*$. Since uniform convergence is preserved under Legendre transform (see e.g. \cite{CSY}), the proof is complete.
\end{proof}

\vspace{2mm}

\begin{ex}\label{Subtlety}
Consider for example the case of a two-dimensional face in $\mathbb{R}^5$ that is a square, such that $\tilde{u}$ takes the value $0$ at three of the vertices and $\epsilon > 0$ at the last. Then $\tilde{u}$ vanishes on the triangle formed by the three vertices where $\tilde{u} = 0$, and on the other triangle agrees with the linear function that vanishes on the long edge and takes value $\epsilon$ at the remaining vertex.
\end{ex}

\begin{ex}\label{SmoothBetween}
Consider the solutions $u_{\epsilon}$ on $\mathbb{R}^3$ to 
$$\det D^2u_{\epsilon} = 1 + \epsilon(\delta_{e_3} + \delta_{-e_3})$$
that are asymptotic to $|x|^2/2$. As $\epsilon$ tends to zero these converge uniformly to $|x|^2/2$ by reasoning similar to that used in the proof of Theorem \ref{Stability1}, thus they cannot be linear on the segment connecting $-e_3$ to $e_3$ for $\epsilon$ small.
\end{ex}

\section{Open Questions}\label{OpenProbs}
In this final section we list several open problems and discuss their significance.

\vspace{2mm}

\begin{enumerate}

\item Theorem \ref{Stability1} shows that in the space of global solutions on $\mathbb{R}^n$ to
$$\det D^2u = 1 + \sum_{i = 1}^M a_i\delta_{p_i}$$
that are asymptotic to $|x|^2/2$, which can be identified with points on an explicit orbifold parametrized by $a_i$ and $p_i$ (\cite{JX}), the set of ``maximally singular" solutions is not small (it has nonempty interior). It is thus natural to ask about the boundary of this set. In particular, are there sharp algebro-geometric conditions on the masses $a_i$ and their locations $p_i$ that guarantee the absence of singularities?

\vspace{2mm}

\item What are the asymptotics of $D^2u$ near the vertices of $P$ in our examples? In the case of a single point mass, the smoothness of the tangent cone to the graph of $u$ away from its vertex was recently established in \cite{HLX} (the solution is not required to be global for this result). A reasonable first goal would be to consider the case of two point masses ($P$ is a line segment) and axisymmetry in $\mathbb{R}^3$, and to study the regularity of the tangent cone to $u$ at a mass. In particular, is the tangent cone smooth away from a ray?

\vspace{2mm}

\item The approach of generating singular Monge-Amp\`{e}re metrics by solving an obstacle problem is quite flexible, and may give a useful perspective on metrics that arise in the large complex structure limit in the study of the Strominger-Yau-Zaslow conjecture. For example, choosing an obstacle that is quadratic when restricted to three rays from the origin and infinity otherwise seems to yield metrics with a singular structure similar to that appearing in \cite{LYZ}. It would be interesting to clarify this connection, and to find other singular structures that can be obtained with our approach, with an eye towards developing intuition for SYZ.

\vspace{2mm}

\item Our examples can also be viewed as solutions to certain geometric optics problems. Generalized versions of such problems correspond to more complicated Monge-Amp\`{e}re type equations known as prescribed Jacobian equations (see e.g. \cite{T}, \cite{JT}, \cite{LT}). It could be interesting to find analogues of our examples for such problems.

\end{enumerate}



\end{document}